\documentclass[11pt]{article}
\usepackage{cite}
\usepackage[cmex10]{amsmath}
\usepackage{algorithm}
\usepackage{algorithmic}
\usepackage{array}
\usepackage{fullpage}
\usepackage{epsfig}
\usepackage{bbm}
\usepackage{enumerate}
\usepackage{amssymb}
\usepackage{amsfonts}
\usepackage{subfig}
\usepackage{multirow}
\usepackage{color}
\usepackage{braket}
\usepackage{url}

\newcommand{\norm}[1]{\left\lVert\,#1\,\right\rVert}

\newcommand{\abs}[1]{\left\vert #1 \right\vert}

\newcommand{\R}{\mathbbm{R}}
\newcommand{\Z}{\mathbbm{Z}}

\newcommand{\minimize}{\text{minimize}}

\newcommand{\argmin}{\text{argmin}}
\newcommand{\st}{\text{subject to}}

\newtheorem{definition}{Definition}

\title{A Fair Assignment of Drivers to Parking Lots}

\author{Nicole Taheri%
	 \thanks{IBM Research Ireland, Damastown Industrial Estate, Dublin 15,
		Ireland. Email: {\tt nicole.taheri@ie.ibm.com}.} 
	 \and Jia Yuan Yu\footnotemark[1]
	 \and Robert Shorten\footnotemark[1] 
}

\begin{document}

\maketitle

\begin{abstract}
	Searching for a parking spot can waste time and gasoline. This waste can be
	reduced by assigning drivers to parking lots based on their destination and
	arrival time. In such a system, drivers could request a parking spot in
	advance and be alerted (e.g., via their phone or vehicle) of their assignment
	to a specific parking lot or available spot. 
%
%
	In this paper, a parking assignment system is described to allocate parking
	spaces in a fair and equitable manner.  Heuristics are developed to solve the
	underlying large scale optimization problem. The efficacy of the system is
	demonstrated by applying our algorithms to real data sets.
\end{abstract}

\let\thefootnote\relax\footnotetext{The authors would like to thank TAPAS
Cologne for letting us use their data.} 
			 
\section{Introduction}
\label{section::intro}

A huge societal problem arises from drivers searching for a parking space.
Statistics illustrating the waste associated with this issue border on the
unbelievable. For example, in the United Kingdom, the average person spends
2,549 hours looking for a parking spot over their lifetime
\cite{telegraph2013}. Further, a survey conducted by Allianz Insurance
estimated that 95\% of British people added an extra mile to their trip
searching for a parking spot \cite{allianz2011}. Such statistics arise
elsewhere and are not just in the UK. For example, it was recently reported
that over one year in a small Los Angeles business district, cars cruising for
parking burned 47,000 gallons of gasoline and produced 730 tons of carbon
dioxide \cite{shoup2007}. Meanwhile, the consulting firm McKinsey recently
claimed that the average car owner in Paris spends four years of his or her
life searching for a parking space \cite{mckinsey}. The parking assignment
problem associated with electric vehicles (EVs) becomes even more acute. Due to
the limited range of these vehicles, the marginal cost of expending energy to
search for spaces may, in some cities, be prohibitively high. Thus, there is a
real and compelling societal and economic need to revisit parking.

Currently, parking guidance systems can alert drivers to the availability of
parking spots in given areas \cite{parkingguidance}. However, this type of
system may not work in a crowded area with few available spots. Moreover, if
parking spaces are not specifically assigned to vehicles, then a spot may no
longer be available by the time the vehicle arrives. Assigning drivers to
parking lots could reduce the time and gasoline wasted by drivers searching for
parking spots and could ensure that the allocation of drivers to parking spots
meets specified criteria.

An assignment of drivers to parking lots can be seen as a resource allocation
problem, where the resource to be distributed is parking spots. An important
quality of resource allocation is fairness; the assignment of goods to users
should be \emph{fair} given a specified definition of fairness \cite{lan2010}.
For a parking lot assignment, fairness can be measured in terms of the travel
time between user destinations and assigned parking lots.

Following on from our prior work in this area \cite{schlote2014,griggs2015}, we
propose an optimization-based method to assign vehicles to parking spots in a
fair way. The principal conceptual difference between our previous work and the
work presented here is that we consider a centralized system that solves an
exact optimization problem (with some relaxation), and is built on an elaborate
notion of fairness from the viewpoint of the vehicle owner. Specifically, we
define a driver-centric measure to determine the fairness of a parking lot
allocation, and construct an algorithm to optimize this fairness measure. We
then test our method on two sets of real driving data: one from Cologne,
Germany \cite{colognedata1,colognedata2} and one from New York City taxis
\cite{nyc-taxi1, nyc-taxi2}. Our results show that using our method could
improve the fairness of a parking assignment by around 140\% when compared to
basic methods with no optimization.

\subsection{Related Work}

Within the research community, questions concerning how to manage parking space
supply-demand mismatch are actively being investigated from a variety of
different angles. One important aspect concerns the delivery of up-to-date,
accurate, real-time information to parking systems in order to achieve greatest
system efficiency. Parking guidance and information (PGI) systems can help
inform drivers of where there are available spots; for example, a system could
list the number of available spots on each floor of a parking garage or in each
parking lot in a city \cite{parkingguidance}. Such a system requires detecting
vehicles in parking spots, and there have been a number of works that consider
the technology (i.e., hardware and/or software) required to detect vehicles in
parking spaces, including \cite{wang2011,benson2006,hodel2003,lee2008}. The
effects of parking guidance on a city have also been analyzed
\cite{benenson2008,rodier2010}. However, these systems do not specifically
assign vehicles to parking spots, they only make drivers aware of the
availability in the specified areas. This type of system is not ideal for
street parking, or situations when there are only a few parking spaces
available, as the availability of a parking space in a crowded area can change
quickly. This type of system could also result in unnecessary overcrowding of
certain areas

A predictive approach is proposed by Caliskan et al.\ \cite{caliskan2007} and
further studied in \cite{klappenecker2010}. The authors develop a method to
predict the likelihood of a parking space being available at the estimated time
that the car will arrive there. In \cite{geng2013}, Geng et al.\ view the
parking problem as a dynamic resource-allocation problem. In this paper,
algorithms are developed to assign drivers to parking lots dynamically as they
leave their destinations. The objective of their formulation is to minimize the
sum of the utility functions for all the drivers. Chou et al. \cite{chou2008}
use a different model and consider a network of parking lots with negotiable
prices; parking lots are selected for drivers in a way that benefits both the
drivers and the car park operators.  The work of Teodorovic and Lucic
\cite{teodorovic2006} proposes a method to define rules for assigning drivers
to parking spots using fuzzy logic and integer programming techniques. The
parking assignment problem has also been solved using a stable matching
algorithm \cite{ayala2012a}, which yields a Nash equilibrium solution, and
through differential pricing \cite{ayala2012b}.  Parking lots can also be
assigned to drivers through matching under preferences \cite{manlove2013}.
Griggs et al. created a distributed privacy-preserving scheme to allocate
parking efficiently with quality of service guarantees \cite{griggs2015}, where
fairness is one such objective. Finally, in \cite{schlote2014} cars are
assigned to lots in a balanced manner, where fairness is considered from the
view point of the parking lot owner.  Our work follows the prior work described
in \cite{griggs2015,schlote2014}. However, there are important differences.  We
provide an advanced driver-centric notion of fairness, and give a heuristic
that solves a very complicated optimization problem that goes far beyond that
described in any of the aforementioned references. Importantly, our fairness
notion can be extended to include other aspects of fairness (aspects associated
with EV ownership) as a particular need arises.

\section{Problem Statement}

We propose an optimization-based method to assign parking lots to drivers in a
fair way.  We assume that the starting and ending times of all trips are known
at the beginning of the day (i.e., a 24-hour period), as well as the starting
and ending locations. 

Without loss of generality, we assume that a vehicle can be collected from a
parking lot as well as dropped off (i.e., a driver may pick up their vehicle
from a specified parking lot and park it at a different parking lot near
his/her destination). We assume that a driver will pick up a vehicle from the
parking lot closest to their origin location and our method will assign the
driver to a parking lot in which to drop off the vehicle that is near their
destination.

\subsection{Notation}
\label{section::notation}

Table \ref{table::notation} defines the terminology to be used in the remainder
of the work.

\begin{table}[h]
\begin{center}
\caption{Sets, Constants and Variables}
\label{table::notation}
\begin{tabular}{ll|l}
	\multicolumn{3}{c}{Sets}  \\
	\hline
  Name            &&  Description \\
	\hline            
  $\mathcal{R}$   && Drivers     \\  
  $\mathcal{T}$   && Time Periods  \\
  $\mathcal{L}$   && Parking Lots  \\ 
	\hline \\
	\multicolumn{3}{c}{Constants}  \\
	\hline
  \multicolumn{2}{c|}{Name} &  Description \\
	\hline
  $M$       & $\in \R^{2\times\abs{\mathcal{L}}}$ & Location of Parking Lots                      \\
  $\bar{x}$ & $\in \Z^{\abs{\mathcal{L}}}$        & Maximum Number of Parking Spots in each Lot        \\
  $s,d$     & $\in \R^{2\times\abs{\mathcal{R}}}$ & Start and End Location of Driver Trips             \\
  $t^s,t^d$ & $\in \R_+^{\abs{\mathcal{R}}}$      & Start and End Time of Driver Trips                 \\
  $\alpha$  & $\in \R+$                           & Average Walking Speed (in miles per hour) \\
	\hline \\
	\multicolumn{3}{c}{Variables}  \\
	\hline
  \multicolumn{2}{c|}{Name} &  Description \\
	\hline
  $y$      & $\in \{0,1\}^{\abs{\mathcal{L}}\times \abs{\mathcal{R}}}$ & Destination Lot Assignment                                    \\  
  $x(t)$   & $\in \Z^{\abs{\mathcal{L}}}$                              & Number of Filled Parking Spots at time $t$                    \\
  $Y(t)$   & $\in \Z^{\abs{\mathcal{L}}}$                               & Number of drivers arriving at each lot at time $t$            \\  
  $Z(t)$   & $\in \Z^{\abs{\mathcal{L}}}$                               & Number of drivers leaving each lot at time $t$                \\  
  $D$      & $\in \R_+^{\abs{\mathcal{R}}}$                            & Distance from Assigned Parking Lot to Destination (in miles) \\
  $\beta$  & $\in \R_+^{\abs{\mathcal{R}}}$                            & Total Overhead Time for each Driver Trip (in hours)                  
\end{tabular}
\end{center}
\end{table}

For all trips $r\in\mathcal{R}$ and parking lots $\ell\in\mathcal{L}$, we
define the binary decision variable $y \in
\{0,1\}^{\mathcal{L}\times\mathcal{R}}$ and the binary parameter $z \in
\{0,1\}^{\mathcal{L}\times\mathcal{R}}$. The parameter $z_{\ell r}$ will
specify from which parking lot the vehicle will depart (its origin), and the
variable $y_{\ell r}$ will specify in which parking lot the driver will
park (its destination).  Because we assume that each vehicle will start in the
parking lot closest to the driver's origin location, $z$ is a parameter, not a
variable.  However, the lot in which a vehicle will park near the destination
is a decision variable to be optimized. The values of $y$ and $z$ can be
defined with the equations:
\begin{equation*}
y_{\ell r} = \begin{cases}
	1 & \text{if trip $r$ ends up in lot $\ell$} \\
	0 & \text{otherwise} 
\end{cases}, \qquad
z_{\ell r} = \begin{cases}
	1 & \text{if trip $r$ starts from lot $\ell$} \\
	0 & \text{otherwise}. 
\end{cases} 
\end{equation*}
The values of $z$ can easily be calculated given the origin locations of
the trips, $s$.

We also define the state variables $Y,Z \in \Z^{\abs{\mathcal{L}}\times
\abs{\mathcal{T}}}$ that keep track of the number of vehicles arriving and
leaving each parking lot over time:
\begin{align*}
Y_{\ell}(t) = \sum_{r \in \mathcal{R}: t^d_r=t} y_{\ell r}, \qquad Z_{\ell}(t) = \sum_{r \in \mathcal{R}: t^s_r=t} z_{\ell r},
\end{align*}
where $t_r^s, t_r^d$ are the origin and destination locations of driver $r$
respectively.  The variable $x_{\ell}(t)$ represents the number of filled (or
parked-in) spots in lot $\ell$ at time $t$, and $\bar{x}_{\ell}$ is the maximum
number of available spots in parking lot $\ell$.  A feasible set of lot
assignments $y_{\ell r}$ and parking lot states $x_{\ell}(t)$ will be in the
set:
\begin{equation}
	\Omega = \left\{ (x,y) : \quad
		\begin{aligned}
	   &x_{\ell}(t) = x_{\ell}(t-1) + Y_{\ell}(t) - Z_{\ell}(t) & \forall \ell\in\mathcal{L}, t\in\mathcal{T}              \\
	   &x_{\ell}(t) \le \bar{x}_{\ell}                          & \forall \ell \in\mathcal{L}, t\in\mathcal{T} \\
	   &\sum_{\ell\in\mathcal{L}} y_{\ell r} = 1                & \forall r\in\mathcal{R}   \\
	   & y_{\ell r}    \in \{0,1\}                              & \forall \ell\in\mathcal{L}, r\in\mathcal{R} 
		\end{aligned}
	\right\}.
	\tag{FEAS}
	\label{eq::feasible-set}
\end{equation}
The first equation in \eqref{eq::feasible-set} updates the state of each
parking lot over time: the number of occupied parking spots is equal to the
number of occupied spots in the previous time period, plus the arriving vehicles
and minus the departing vehicles in the previous period. The second equation in
\eqref{eq::feasible-set} ensures the number of vehicles in each lot is less
than the capacity, the third equation makes sure that each trip ends in only one
parking lot, and the last equation constrains $y$ to be binary.

\section{Optimizing Fairness}
\label{section::optimizing-fairness}

Assigning drivers to parking lots leads us to the following questions: Is there
a \emph{fair} way to assign drivers to parking lots? And what does it mean for
an assignment of parking lots to be fair?  

In our problem, a fair assignment will be one in which the distance between
drivers' destinations and their assigned parking lots are similar. Let $\beta_r$
be the amount of time driver $r \in \mathcal{R}$ spends traveling from the
assigned parking lot to his/her destination:
\begin{equation}
	\beta_r = \alpha \norm{d_r - \sum_{\ell\in\mathcal{L}} (y_{\ell r} \cdot M_{\ell}) }_1 .
\end{equation}
We define fairness as a lack of \emph{envy} among drivers, where the envy
between two drivers  $r_1, r_2 \in \mathcal{R}$ is defined as:
\begin{equation}
	E_{r_1, r_2} := \abs{\beta_{r_1} - \beta_{r_2}},
	\label{eq::driver-envy}
\end{equation}
where $\beta_r$ is the walking time for driver $r \in \mathcal{R}$ from his/her
parking lot to destination. 

We define fairness as the mean value of $E_{r_1, r_2}$ over all pairs of
drivers  $r_1, r_2 \in \mathcal{R}$, where a smaller value is better. To attain
an assignment of drivers to parking lots that optimizes this fairness measure,
we can minimize the mean of the envy values over all pairs of drivers with the
objective function:
\begin{equation}
	F(\beta) = \underset{r_1, r_2 \in \mathcal{R}}{\text{mean}} (E_{r_1,r_2} )
	= (1/\abs{R}^2) \sum_{r_1, r_2 \in \mathcal{R}} \abs{\beta_{r_1} - \beta_{r_2}}.
	\label{eq::driver-fairness}
\end{equation}
\begin{definition}
	A \emph{fair assignment} of destination lots to drivers results in a value of 
	$F(\beta)$ that is as small as possible. 
\end{definition}

\subsection{Formulation}
\label{section::formulation}

We assume that we are given the origin and destination locations and times for
each driver, $(s,d,t^s,t^d)$ at the start of the day (i.e., a 24-hour period).
Let $g: \R^{\abs{\mathcal{R}}} \rightarrow \R$ be a function that takes the
waiting time of each driver $r \in \mathcal{R}$ and assesses the fairness,
where a smaller value implies a more fair assignment. Then we can formulate the
problem that maximizes fairness as:

\begin{align}
	\underset{x,y,\beta}{\minimize} \quad & g(\beta)  \tag{OPT-FAIR} \label{eq::fairness-opt} \\
	      \st \quad &(x,y) \in \Omega \notag \\
									& \beta_r = \alpha \norm{d_r - \sum_{\ell\in\mathcal{L}} (y_{\ell r} \cdot M_{\ell}) }_1 & \forall r\in \mathcal{R} \label{eq::distance} 
\end{align}

The first constraint ensures that the lot assignments defined by $x$ and $y$
are feasible for the given drivers and parking lots.  Constraint
\eqref{eq::distance} calculates the distance from the driver's destination to
the assigned parking lot and multiplies it by the average walking speed,
$\alpha$, to get the total travel time for driver $r$.  The formulation
\eqref{eq::fairness-opt} is a mixed-integer linear program (MILP) that can be
solved in CPLEX \cite{cplex}. 

This optimization problem can be solved for a chosen time period (e.g, a single
day); our aim is to choose an objective function $g(\beta)$ for
\eqref{eq::fairness-opt}, so that the average overhead waiting time over a set
of drivers is fair.

\subsection{Intractability of $F(\beta)$ as the objective}
\label{section::intractability}

Using the definition of fairness based on the function $F(\beta)$ in
\eqref{eq::driver-fairness}, an optimally fair lot assignment can be found by
setting the objective function in \eqref{eq::fairness-opt} to $g(\beta) =
F(\beta)$.  That is, ideally we could minimize the mean of all the envy values, or
differences between driver walking times.  However, setting $g(\beta) =
F(\beta)$ as the objective in \eqref{eq::fairness-opt} results in an
intractable optimization problem; taking the absolute value of each pairwise
difference in travel times results in a large problem that cannot be solved
in reasonable time.  For example, we used CPLEX to try to solve the problem
\eqref{eq::fairness-opt} with the objective $g(\beta) = F(\beta)$, for
$\abs{\mathcal{L}} = 10$ parking lots, $\abs{\mathcal{R}} = 100$ drivers, and
$\abs{\mathcal{T}} = 24$ time periods, and the solver did not even find a
feasible solution after 12 hours. 

\subsection{Choosing an Objective Function}
\label{section::choosing-objective}

Given the intractability of using $F(\beta)$ as the objective, we would like to
construct a method that achieves a similar result, but that is computationally
tractable. The idea of minimizing $F(\beta)$ in \eqref{eq::fairness-opt} is to
obtain a lot assignment in which all drivers' walking times are similar. We
propose a method to find a solution of \eqref{eq::fairness-opt} that achieves
this same goal, but with an objective that makes this problem tractable.

Given a feasible set of lot assignments $y \in \{0,1\}^{\abs{L} \times
\abs{\mathcal{R}}}$ with corresponding walking times $\beta \in
\R^{\abs{\mathcal{R}}}$, the mean walking time is defined as:
\begin{equation}
	H(\beta) := (1/\abs{\mathcal{R}}) \sum_{r \in \mathcal{R}} \beta_r.
	\label{eq::define-mean}
\end{equation}
A fair assignment of drivers to parking lots can be attained by minimizing the
absolute value of the differences between each driver's walking time and the
mean walking time, $H(\beta)$. This will provide a solution where the
walking times of the drivers ares similar and as close to the mean value
as possible. 

Given a feasible solution $(\widehat{x},\widehat{y},\widehat{\beta})$ to
\eqref{eq::fairness-opt}, the mean walking time $H(\widehat{\beta})$ can be calculated
and the problem \eqref{eq::fairness-opt} can be re-solved with the objective: 
\begin{equation}
	G(\beta) = \sum_{r \in \mathcal{R}} \abs{\beta_r - H(\widehat{\beta})}.
	\label{eq::fair-objective}
\end{equation}
Our method is based on solving \eqref{eq::fairness-opt} with this objective,
which will yield a tractable problem and a fair solution. 

To further reduce the size of the problem, we can consider only a subset of the
lot assignments as decision variables. That is, because some of the walking
times $\beta_r$ may be sufficiently close to the mean value
$H(\widehat{\beta})$ already, we can fix the lot assignments within a chosen
range of $H(\widehat{\beta})$ and only change those with a large or small
travel time. Define the set $S\subseteq\mathcal{R}$:
\begin{equation}
	S := \{r \in \mathcal{R}: \widehat{\beta}_r \in [ (1-\epsilon) H(\widehat{\beta}), (1+\epsilon) H(\widehat{\beta}) ] \},
\end{equation}
where $\epsilon > 0$ is a chosen value that determines the desired range of
walking times from the mean. 

For a feasible set of lot assignments, $\widehat{y} \in
\{0,1\}^{\abs{\mathcal{L}} \times \abs{\mathcal{R}}}$ with corresponding
walking times $\widehat{\beta} \in \R^{\abs{\mathcal{R}}}$, we define the set
of constraints based on $S$ as:
\begin{equation}
	\widehat{\Omega}(S) = 
	\Set{
		(x, y, \beta) | 
	\begin{aligned}
		\begin{array}{ll}
			(x, y) \in \Omega  \\
	    \beta_r = \alpha \cdot \norm{d_r - \sum_{\ell\in\mathcal{L}} (y_{\ell r} \cdot M_{\ell}) }_1 & \forall r\in \mathcal{R} \\
			y_{\ell r} = \widehat{y}_{\ell r} & \forall r \in S, \ell \in \mathcal{L} \\
			\beta_{r} = \widehat{\beta}_{r} & \forall r \in S
		\end{array} \\
	\end{aligned}
	},
	\label{eq::define-omega-hat}
\end{equation}
and we define the new objective function:
\begin{equation*}
	G(\beta,S) = \sum_{r \not\in S} \abs{\beta_r - H(\widehat{\beta})}.
\end{equation*}
The following optimization problem will then find new lot assignments that
optimizes fairness and reduces the envies among drivers' walking times:
\begin{align}
	\underset{x,y,\beta}{\minimize} \quad & G(\beta,S) \tag{OPT-MEAN} \label{eq::fairness-mean} \\
	      \st \quad &(x,y,\beta) \in \widehat{\Omega}(S) \notag 
\end{align}

Notice that after solving the MILP \eqref{eq::fairness-mean} and obtaining a
new assignment $(y^{\star},x^{\star},\beta^{\star})$, the mean walking time
will change according to the new lot assignments and walking times
$\beta^{\star}$.  That is, the mean walking of the new solution,
$H(\beta^{\star})$, will differ from $H(\widehat{\beta})$, and thus
\eqref{eq::fairness-mean} can be re-solved with the new mean walking time,
$H(\beta^{\star})$. This process of re-solving can continue until the mean
walking value converges (i.e., stops changing).  Iteratively re-solving
\eqref{eq::fairness-mean} with the updated mean value $H(\beta)$, and the
corresponding updated set $S$ will reduce the envies of all the drivers until
their values cannot be any closer.

Our full method is described in Algorithm \ref{alg::min-envy}. This method
iteratively minimizes the difference between a subset of the drivers' walking
times and the mean walking time of the previous iteration,
$H(\widehat{\beta})$, until the solution has converged. The algorithm takes
three parameters: $\epsilon$ is the desired range from the mean walking time,
$\delta$ is the tolerance used to determine if the algorithm has converged, and
\texttt{maxiter} is the maximum number of iterations.
\begin{algorithm}
	\caption{Minimizing the envy}
	\label{alg::min-envy}
	\begin{algorithmic}[1]
		\STATE Choose values for $\epsilon, \delta, $ \texttt{maxiter}
		\STATE Find a feasible lot assignment $y^0 \in \{0,1\}^{\abs{L} \times
		\abs{\mathcal{R}}}$ with corresponding walking times $\beta^0 \in \R^{\abs{\mathcal{R}}}$
		\STATE $i \leftarrow 1$, \texttt{converged} $\leftarrow 0$
		\WHILE{not converged}
			\STATE $S^i \leftarrow \{r \in \mathcal{R}: \beta^{i-1}_r \in [ (1-\epsilon) H(\beta^{i-1}), (1+\epsilon) H(\beta^{i-1}) ] \}$
			\STATE $(x^i,y^i, \beta^i) \leftarrow \underset{(x,y,\beta) \in \widehat{\Omega}(S^i)}{\argmin} \; G(\beta,S^i)$ 
			\IF{$\abs{H(\beta^i) - H(\beta^{i-1})} < \delta$ {\bf or} $i > $ \texttt{maxiter}}
				\STATE \texttt{converged} $\leftarrow 1$
			\ENDIF 
			\STATE $i \leftarrow (i + 1)$
		\ENDWHILE 
	\end{algorithmic}
\end{algorithm}

\section{Results}
\label{section::results}

We tested Algorithm \ref{alg::min-envy} on two sets of real driving data. A
description of the data and the results are below. We also compared our method
to the parking assignment formulation considered by Geng and Cassandras in
\cite{geng2012, geng2013}.  All of our implementation was done in Matlab
\cite{matlab} and used IBM ILOG CPLEX \cite{cplex} to solve the optimization
problems.

\subsection{Data}
\label{section::data}

Real trip data was collected from New York City taxis in 2013 \cite{nyc-taxi1,
nyc-taxi2}.  From this dataset, we extracted information on the origin and
destination locations and times for over 50,000 trips in January 2013. 

We also used real driving data from Cologne, Germany, collected for the
Simulation of Urban Mobility (SUMO) project \cite{colognedata1,colognedata2}.
This data is known as the TAPAS Cologne data and provides us with driver
behaviors in Cologne for a 2-hour period.

In order to model a system with known parking lot locations, for each data set
we created 10 parking lots in the given area. For each set of drivers and their
known destinations $d \in \R^{2\times\abs{\mathcal{R}}}$, we used $k$-means
clustering \cite{kmeans} on the locations $d$ to find 10 clusters of
destinations; we then set the locations of the 10 parking lots to the centroids
of each of the clusters. 

Let $w \sim \bar{\mathcal{U}}[a,b]$ signify that $w$ is a pseudorandom integer
from a discrete uniform distribution between $a$ and $b$.  For each of the
parking lots, we assigned a capacity $\bar{x}_{\ell}$ and initial occupancy
$x_{\ell}(0)$ of:
\begin{align*}
	\bar{x}_{\ell} \; &\sim\;  \bar{\mathcal{U}}\left[\left(\abs{\mathcal{R}}/\abs{\mathcal{L}}\right) + 1,
		\quad \left(\abs{\mathcal{R}}/\abs{\mathcal{L}}\right) + 2\right] \\
	x_{\ell}(0) \;    &\sim \; \bar{\mathcal{U}}[(\bar{x}_{\ell}/4), \quad (3\cdot\bar{x}_{\ell}/4)].
\end{align*}
In other words, the capacity of each parking lot will be the average number of
parking spots needed in each lot
$\left(\abs{\mathcal{R}}/\abs{\mathcal{L}}\right)$ plus either 1 or 2, and each
lot will initially be between one-quarter and three-quarters full. The average
walking speed $\alpha$ was assumed to be 5 km/hour \cite{walkingspeed}.

\subsection{Fairness Results}
\label{section::data}


To determine the relative fairness of our method, we also implemented the two
methods described in the algorithms below.  Algorithm \ref{alg::min-sum}
(Minimum Sum) solves the optimization problem \eqref{eq::fairness-opt} once
with the objective of minimizing the sum of the travel times between the
parking lot and destination.  Algorithm \ref{alg::no-scheme} (No Scheme)
simulates what might happen with no scheme, if drivers were to park in the
closest available lot to their destination. (Note that this algorithm will
perform slightly better than no scheme, because it assumes drivers know exactly
where the best available spot is.)

\begin{algorithm}[H]
	\caption{Minimizing the sum}
	\label{alg::min-sum}
	\begin{algorithmic}[1]
		\STATE Find the lot assignments by solving the optimization problem:
		\begin{align*}
		(x,y, \beta) \leftarrow & \underset{(x,y) \in \Omega}{\argmin} \qquad \sum_{r \in \mathcal{R}} \beta_r\\
		& \st  \quad \beta_r = \alpha \cdot \norm{d_r - \sum_{\ell\in\mathcal{L}} (y_{\ell r} \cdot M_{\ell}) }_1 & \forall r\in \mathcal{R} 
		\end{align*}
	\end{algorithmic}
\end{algorithm}

\begin{algorithm}[H]
	\caption{No scheme}
	\label{alg::no-scheme}
	\begin{algorithmic}[1]
		\FOR{$t =1,\dots,\abs{T}$}
			\FOR{all drivers $r \in \mathcal{R}$ arriving at their destination at time $t$}
				\STATE Find the closest lot to the destination of driver $r$ that has
					spaces available and assign the driver to this parking lot
			\ENDFOR 
		\ENDFOR 
	\end{algorithmic}
\end{algorithm}

The figures below show a comparison of the fairness of our algorithm to the
fairness of Algorithms \ref{alg::min-sum} and \ref{alg::no-scheme}. For the
data for New York City, we ran 500 tests on the data, where each test randomly
extracted a different set of 100 drivers over a 24-hour period from the data
set described in Section \ref{section::data}.  For each run, we found the mean
envy value for that run (i.e., we calculated the value of $F(\beta)$ in
\eqref{eq::driver-fairness} for the solution to each method).

\begin{figure}[H]
	\begin{center}
   \subfloat[Progression of $F(\beta)$ over iterations of Algorithm
		\ref{alg::min-envy}]{\label{fig::nyc-results1}
			\includegraphics[scale=0.43]{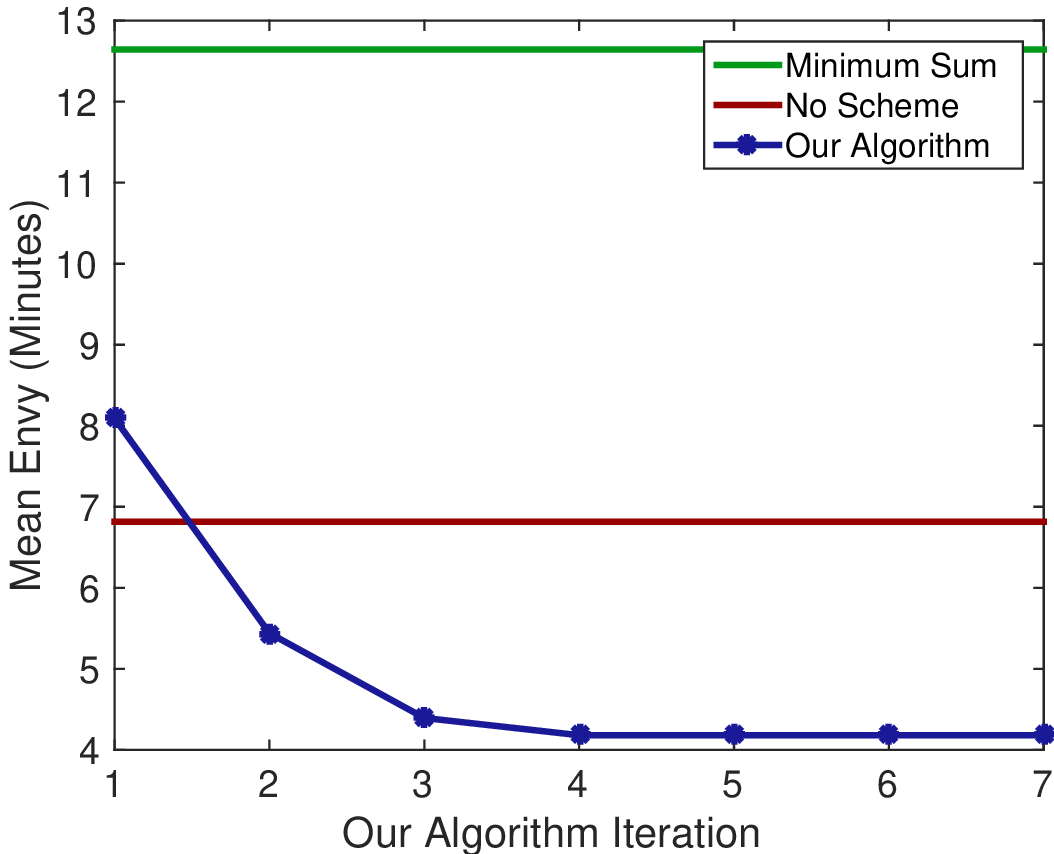}}
		\hspace{0.1in}
    \subfloat[Probability of $F(\beta)$ values]{\label{fig::nyc-results2}
			\includegraphics[scale=0.43]{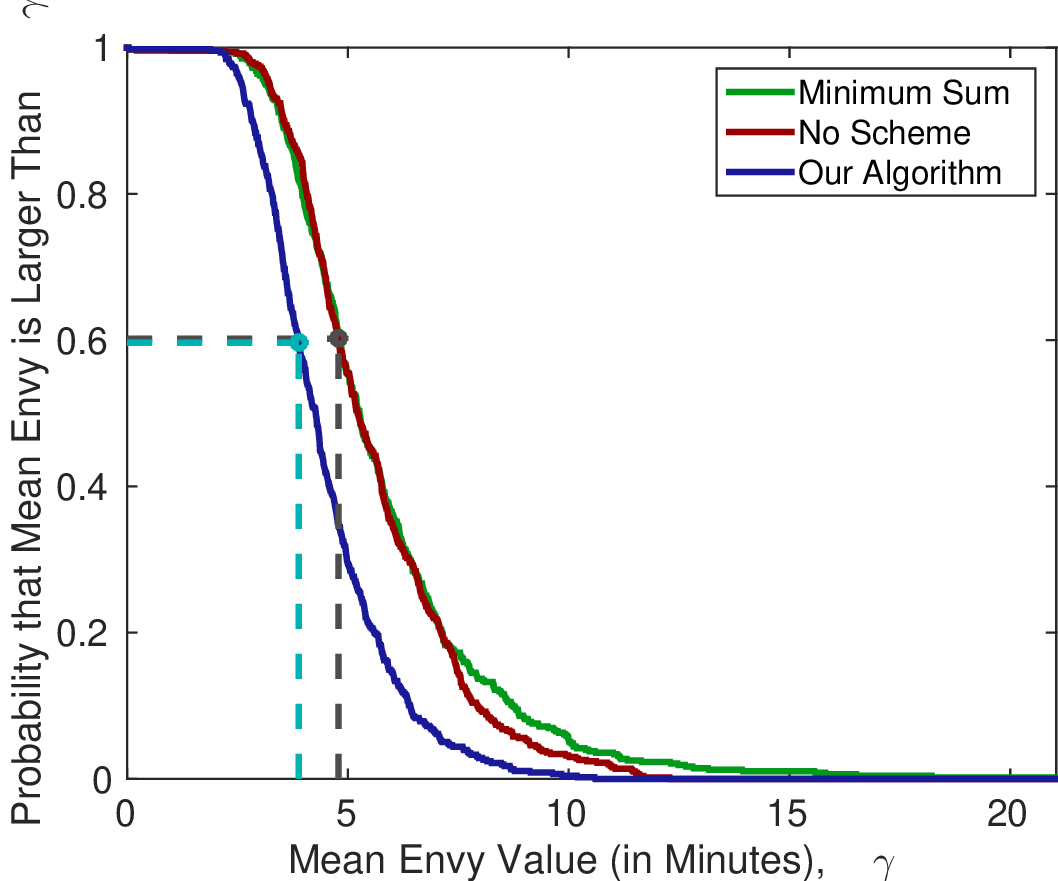}}
	\end{center}
	\begin{center}
     \caption[NYC Taxi Data]
		 {Implementation Results on NYC Taxi Data}
     \label{fig::all-nyc-results}
	\end{center}
\end{figure}

Figure \ref{fig::nyc-results1} shows a progression of the mean envy values,
$F(\beta)$, over iterations of Algorithm \ref{alg::min-envy} in comparison to
Algorithm \ref{alg::min-sum} (Minimum Sum) and Algorithm \ref{alg::no-scheme}
(No Scheme).  Figure \ref{fig::nyc-results2} shows the probability that the
mean envy $F(\beta)$ will be greater than a chosen value. That is, the $x$-axis
in Figure \ref{fig::nyc-results2} corresponds to a given mean envy value
$\gamma$, and the $y$-axis is the probability that the mean envy $F(\beta)$ is
greater than $\gamma$ for each method.  For example, Figure
\ref{fig::nyc-results2} shows that 60\% of envy values for Algorithm
\ref{alg::no-scheme} (No Scheme) are greater than 4.79 minutes, and 60\% of envy
values for Algorithm \ref{alg::min-envy} (Minimum Envy) are greater than 3.89
minutes.  Moreover, the mean improvement of Algorithm \ref{alg::min-envy} over
Algorithm \ref{alg::min-sum} is 28.4\%, and the mean improvement of Algorithm
\ref{alg::min-envy} over Algorithm \ref{alg::no-scheme} is 25.6\%.

The same tests were run on the data from Cologne, described in Section
\ref{section::data}.  We ran 500 tests on the data, where each test used
randomly extracted driving data for 100 drivers over a 2-hour period and the
data was extrapolated to cover a 24-hour period.  

\begin{figure}[H]
	\begin{center}
   \subfloat[Progression of $F(\beta)$ over iterations of Algorithm
		\ref{alg::min-envy}]{\label{fig::cologne-results1}
			\includegraphics[scale=0.43]{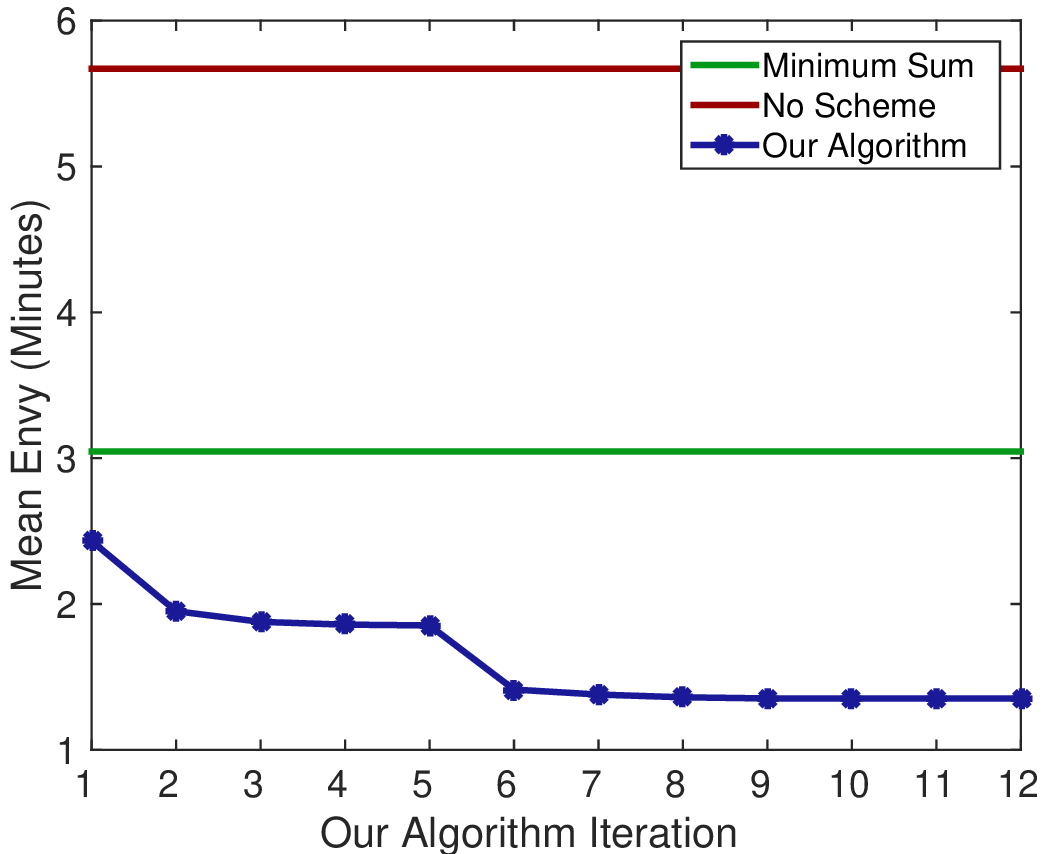}}
		\hspace{0.1in}
    \subfloat[Probability of $F(\beta)$ values]{\label{fig::cologne-results2}
			\includegraphics[scale=0.43]{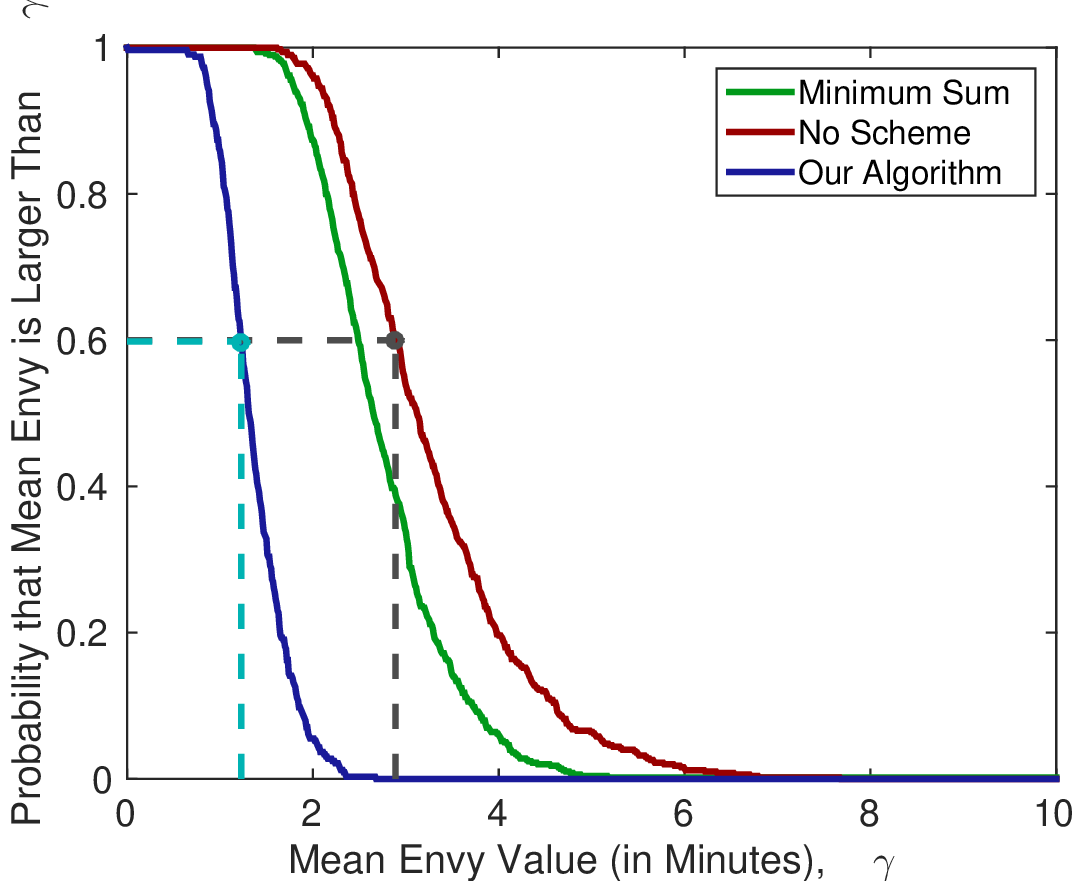}}
	\end{center}
	\begin{center}
     \caption[Cologne Data]
		 {Implementation Results on Cologne Driving Data}
     \label{fig::all-cologne-results}
	\end{center}
\end{figure}

Figures \ref{fig::cologne-results1} and \ref{fig::cologne-results2} show the
results for the Cologne data that are analogous to Figures
\ref{fig::nyc-results1} and \ref{fig::nyc-results2} for the New York City data.
Figure \ref{fig::cologne-results2} shows that 60\% of envy values for Algorithm
\ref{alg::no-scheme} (No Scheme) are greater than 2.89 minutes, and 60\% of envy
values for Algorithm \ref{alg::min-envy} (Minimum Envy) are greater than 1.23
minutes.  Moreover, the mean improvement of Algorithm \ref{alg::min-envy} over
Algorithm \ref{alg::min-sum} is 139.6\%, and the mean improvement of Algorithm
\ref{alg::min-envy} over Algorithm \ref{alg::no-scheme} is 100.2\%.

The Jain's fairness measure is commonly used in congestion control protocols to
measure fairness of a resource allocation \cite{mo2000,lan2010,jainswiki}.  The
value of the Jain's measure is between 0 and 1, where a more fair assignment
results in a higher value.  For a set of travel times $\beta_r$, the Jain's
fairness measure is defined as:
\begin{equation}
	\mathcal{J}(\beta) = \frac{ \left(\sum_{r \in \mathcal{R}} \beta_r \right)^2}{\abs{\mathcal{R}} \cdot \sum_{r \in \mathcal{R}} \beta_r^2 }
	\label{eq::jains}
\end{equation}

\begin{figure}[H]
	\begin{center}
    \subfloat[Jain's measures using NYC Taxi Data]{\label{fig::nyc-jains}
			\includegraphics[scale=0.43]{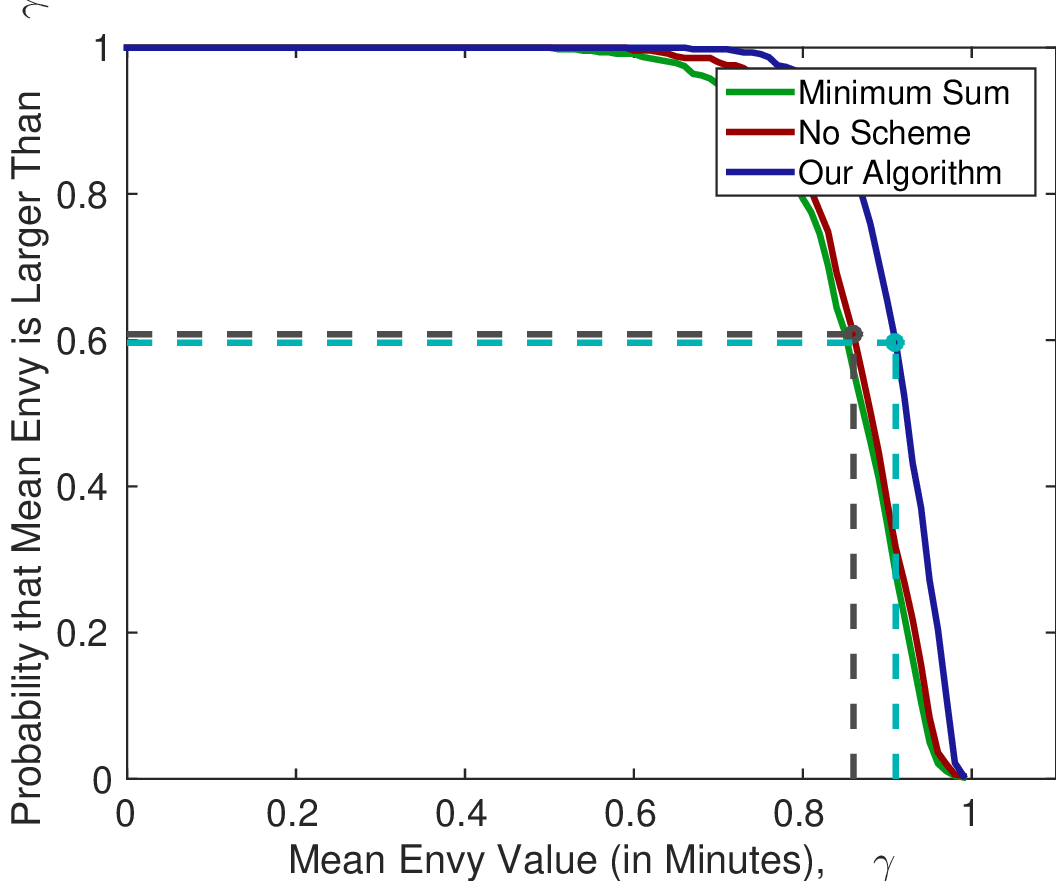}}
		\hspace{0.1in}
    \subfloat[Jain's measures using Cologne Driving Data]{\label{fig::cologne-jains}
			\includegraphics[scale=0.43]{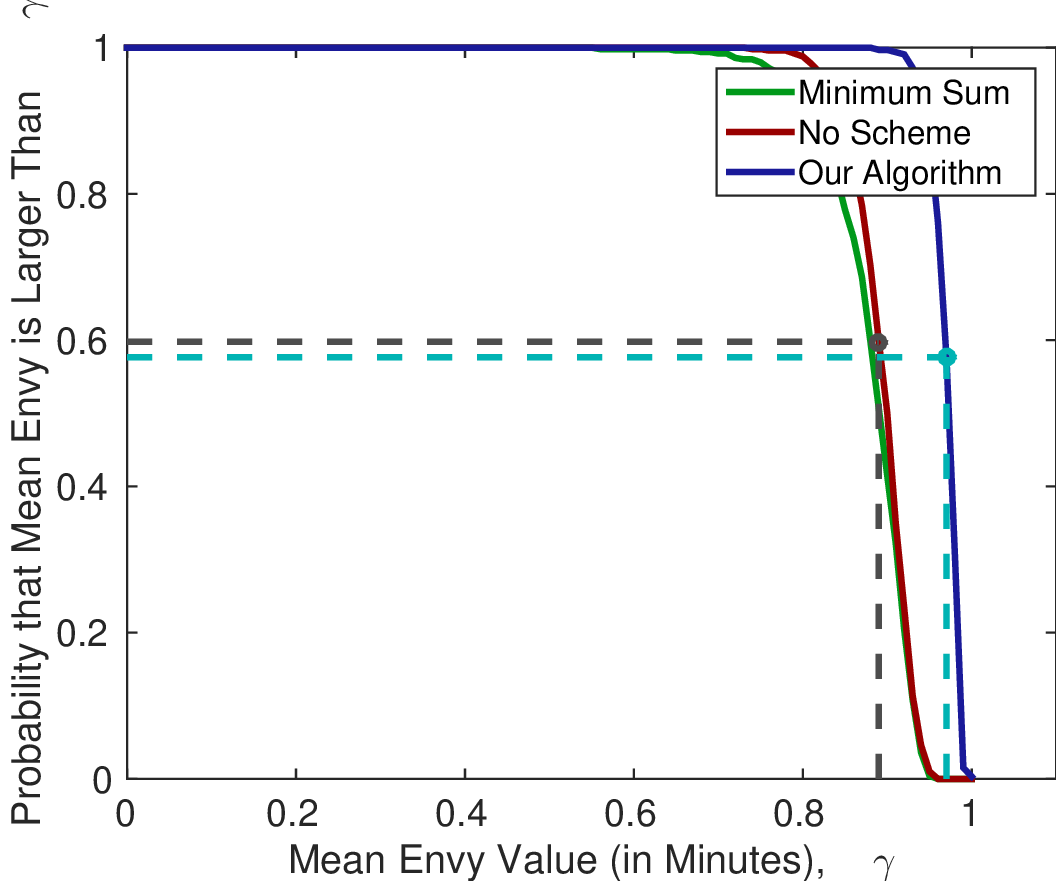}}
	\end{center}
	\begin{center}
     \caption[Jain's Measures]
		 {Probability of Jain's Measure, $\mathcal{J}(\beta)$}
		 \label{fig::all-jains}
	\end{center}
\end{figure}

We computed the Jain's fairness measure for each dataset for each method;
Figure \ref{fig::all-jains} shows the probability that the Jain's measure
$\mathcal{J}(\beta)$ will be greater than a chosen value. 

Note that, as opposed to in Figures \eqref{fig::cologne-results2} and
\eqref{fig::nyc-results2}, for the Jain's fairness measure, it is better for
the values to be higher. 

For Figure \ref{fig::all-jains} the $x$-axis corresponds to a given value
$\gamma$, and the $y$-axis is the probability that the Jain's measure
$\mathcal{J}(\beta)$ is greater than $\gamma$ for each method.  For example,
using the Cologne data, Figure \ref{fig::nyc-results2} shows that 60\% of
Jain's measure values for Algorithm \ref{alg::no-scheme} (No Scheme) are
greater than 0.89, and 60\% of Jain's measure values for Algorithm
\ref{alg::min-envy} (Minimum Envy) are greater than 0.97.  Using the data from
NYC taxi trips, 60\% of Jain's measure values for Algorithm
\ref{alg::no-scheme} (No Scheme) are greater than 0.86, and 60\% of Jain's
measure values for Algorithm \ref{alg::min-envy} (Minimum Envy) are greater
than 0.91.  

Moreover, using the Cologne data, the mean improvement of Algorithm
\ref{alg::min-envy} over Algorithm \ref{alg::min-sum} is 8.95\%, and the mean
improvement of Algorithm \ref{alg::min-envy} over Algorithm
\ref{alg::no-scheme} is 7.67\%.  Using the NYC taxi trips, the mean improvement
of Algorithm \ref{alg::min-envy} over Algorithm \ref{alg::min-sum} is 6\%, and
the mean improvement of Algorithm \ref{alg::min-envy} over Algorithm
\ref{alg::no-scheme} is 4.61\%.

%

\subsection{Comparison to Other Methods}
\label{section::data}

We compared our algorithm to the method used by Geng and Cassandras in
\cite{geng2012, geng2013}. In this work, the authors assign parking lots to
drivers dynamically as they leave their destinations (that is, the driving
patterns are not known ahead of time). Geng et al.\ do not consider fairness of
lot assignments; their method minimizes the sum of the utility functions of
each driver requesting a parking spot at a given time.

We implemented the method of Geng et al.\ as described in the paper
\cite{geng2013} and compared their method to ours. We generated random data
using the same methods they described, according to the specified distributions
and parameters meant to simulate parking around Boston University.  We also
implemented our method on a dynamic basis, so that our algorithm was
implemented to assign parking lots to the drivers that were currently
requesting parking spots.

Geng et al.\ define the utility function of each driver $r \in \mathcal{R}$ as
\begin{equation}
	J_{\ell r} = \lambda_r \frac{M_{\ell r}}{M_r} + (1-\lambda_r) \frac{D_{\ell r}}{D_r},
	\label{eq::define-j-cost}
\end{equation}
where $\lambda_r \in [0,1]$ is a weight defined by the user, $M_{\ell r}$ is
the monetary cost for driver $r$ to park in parking lot $\ell$, and $D_{\ell
r}$ is the distance from driver $r$'s destination to parking lot $\ell$. The
values $M_r$ and $D_r$ are also provided by the user, and these are the maximum
values of $M_{\ell r}$ and $D_{\ell r}$ that user $r$ is willing to accept.

For each time step $k = 1, 2, \dots$, let $P(k) \subseteq \mathcal{R}$ be the
set of drivers that have requested a parking spot after time $k$ and have not
yet arrived at their destination by time $k$.  Using notation consistent with
Table \ref{table::notation} and adding the time step parameter $k$, Geng et
al.\ find a lot assignment $y(k) \in
\{0,1\}^{\abs{\mathcal{L}}\times\abs{\mathcal{P}(k)}}$ by solving the problem
\begin{align}
	\underset{x,y(k)}{\minimize} \quad & \sum_{\ell\in\mathcal{L}, r \in \mathcal{P}(k)} y_{\ell r}(k) \cdot J_{\ell r}  +
		\sum_{r \in \mathcal{P}(k)} \left( 1 - \sum_{\ell \in \mathcal{L}} y_{\ell r}(k) \right) \label{eq::cassandras-opt} \tag{SMART-PARK} \\
	      \st \quad & x_{\ell}(t) = x_{\ell}(t-1) + Y_{\ell}(t) - Z_{\ell}(t) & \forall \ell\in\mathcal{L}, t\in\mathcal{T} \notag \\
	   							& x_{\ell}(t) \le \bar{x}_{\ell}                          & \forall \ell \in\mathcal{L}, t\in\mathcal{T}\notag \\
	   							& \sum_{\ell\in\mathcal{L}} y_{\ell r}(k) \le 1           & \forall r\in\mathcal{P}(k)                  \label{eq::lessthan}      \\
	   							& y_{\ell r}    \in \{0,1\}                               & \forall \ell\in\mathcal{L}, r\in\mathcal{P}(k) \notag \\
									& \sum_{\ell \in \mathcal{L}} (y_{\ell r}(k) \cdot J_{\ell r}) \le 
									  \sum_{\ell \in \mathcal{L}} (y_{\ell r}(k-1) \cdot J_{\ell r}) \label{eq::cost-improvement}
\end{align}
Notice that the only difference between the first four constraints in
\eqref{eq::cassandras-opt} and the set defined by $\Omega$ in
\eqref{eq::feasible-set} is the equation \eqref{eq::lessthan}, where this
equation is an inequality in \eqref{eq::cassandras-opt} and an equality in
\eqref{eq::feasible-set}. This implies that each driver may not necessarily be
assigned to a parking spot in \eqref{eq::cassandras-opt}, but is guaranteed a
parking spot assignment with the constraints in \eqref{eq::feasible-set}.  The
last constraint \eqref{eq::cost-improvement} ensures that any new assignment to
driver $r$ is no worse than their assignment from the previous period, $k-1$.

To compare our method, we used the formulation \eqref{eq::cassandras-opt}, but
changed the objective function to:
\begin{align}
	\underset{x,y(k)}{\minimize} \quad & \sum_{\ell\in\mathcal{L}, r \in \mathcal{P}(k)} \left(y_{\ell r}(k) \cdot J_{\ell r}  - \bar{J}(k)\right) +
		\sum_{r \in \mathcal{P}(k)} \left( 1 - \sum_{\ell \in \mathcal{L}} y_{\ell r}(k) \right),
	\label{eq::cass-fair-obj}
\end{align}
where 
\begin{equation*}
	\bar{J}(k) :=  \frac{1}{\abs{\mathcal{P}(k)}} \left( \sum_{\ell\in\mathcal{L}, r \in \mathcal{P}(k)} y_{\ell r}(k-1) \cdot J_{\ell r} \right).
\end{equation*}
This is analogous to the objective defined in \eqref{eq::fair-objective},
where the only difference is that \eqref{eq::cass-fair-obj} substitutes
the travel times $\beta$ with the utility function values $J$.  We then use the
iterative method defined in Algorithm \ref{alg::min-envy} with the optimization
problem defined by the objective \eqref{eq::cass-fair-obj} and the constraints
to \eqref{eq::cassandras-opt}. This allows us to compare the fairness of our
algorithm with that of Geng et al. 

Table \ref{table::cass-compare} compares the relative improvement in the mean
of the envy values $F(\beta)$ as defined in \eqref{eq::driver-fairness}, over
500 runs.  Each run uses different randomly generated data using the guidelines
provided in \cite{geng2013}, and compares the fairness of Geng et al.'s method
to Algorithm \ref{alg::min-envy} using the formulation
\eqref{eq::cassandras-opt} with the objective \eqref{eq::cass-fair-obj}. The
values in Table \ref{table::cass-compare} are the mean of the relative
improvement in the fairness measures $F(\beta)$ and $\mathcal{J}(\beta)$ for
the solutions of each of the 500 runs.

\begin{table}[h]
	\begin{center}
	\begin{tabular}{|c|c|}
		\hline
		Mean Envy & Jain's Fairness Measure \\
		\hline
		14.51\% &  4.56\% \\ 
		\hline
	\end{tabular}
	\caption{Mean Relative Improvement in Fairness using Algorithm \ref{alg::min-envy} 
	over \eqref{eq::cassandras-opt}}
	\label{table::cass-compare}
	\end{center}
\end{table}

\section{Conclusion}
\label{section::conclusion}

In this work, we constructed a measure to determine the fairness of a given
assignment of drivers to parking lots. We used this measure to construct an
optimization-based algorithm that iteratively reduces the envy among drivers
(or iteratively increases fairness) and finds a fair parking lot allocation.
Our method was tested on real driving data and compared to a number of similar
methods to show that our algorithm provides fair parking lot assignments.

\nocite{*}
\bibliographystyle{unsrt}
\bibliography{parking_refs}

\end{document}